\def\d{{\delta}}

\def\t{{\tau}}
\def\Th{{\Theta}}
\def\w{{\omega}}
\def\W{{\Omega}}

\def\bB{{\mathbf B}}
\def\bG{{\mathbf G}}

\def\E{{\mathbf E}}

\def\G{{\cal G}}

\def\K{{\cal K}}

\def\mt{{\emptyset}}
\def\inf{{\infty}}
\def\ra{{\rightarrow}}

\def\l({{\left(}}
\def\r){{\right)}}
\def\cpf{{$\diamond$}}
\def\pf{{$\Box$}}
\def\"{{^\prime}}

\def\lc{{\lceil}}
\def\rc{{\rceil}}
\def\lf{{\lfloor}}
\def\rf{{\rfloor}}
\def\({{\Biggl(}}
\def\){{\Biggr)}}
\def\[{{\Biggl[}}
\def\]{{\Biggr]}}

\documentclass[12pt]{article}
\newtheorem{theorem}{Theorem}

\usepackage{latexsym,psfig,multicol,amssymb}

\begin{document}

%
%
\title{Pebbling in Dense Graphs}
\author{
Andrzej Czygrinow~\\
Department of Mathematics and Statistics\\
Arizona State University\\
Tempe, Arizona 85287-1804\\
email: andrzej@math.la.asu.edu\\
 \\
and \\
 \mbox{}\\
Glenn Hurlbert\thanks{Partially supported by National Security Agency
        grant \#MDA9040210095.}~\\
Department of Mathematics and Statistics\\
Arizona State University\\
Tempe, Arizona 85287-1804\\
email: hurlbert@asu.edu\\
}
\maketitle
\newpage

%
%
\begin{abstract}
A configuration of pebbles on the vertices of a graph is solvable if one
can place a pebble on any given root vertex via a sequence of pebbling steps.
The pebbling number of a graph $G$ is the minimum number $\pi(G)$ so that
every configuration of $\pi(G)$ pebbles is solvable.
A graph is Class 0 if its pebbling number equals its number of vertices.
A function is a pebbling threshold for a sequence of graphs if a randomly
chosen configuration of asymptotically more pebbles is almost surely
solvable, while one of asymptotically fewer pebbles is almost surely not.
Here we prove that graphs on $n\ge 9$ vertices having minimum degree at least
$\lf n/2\rf$ are Class 0, as are bipartite graphs with $m\ge 336$ vertices in
each part having minimum degree at least $\lf m/2\rf +1$.
Both bounds are best possible.
In addition, we prove that the pebbling threshold of graphs with minimum
degree $\d$, with $\sqrt{n} \ll \d$, is $O(n^{3/2}/\d)$, which is tight when $\d$ is proportional
to $n$.
\vspace{0.2 in}

\vspace{0.2 in}

\end{abstract}

\newpage

%
%
\section{Introduction}\label{Intro}

%
%
\subsection{Pebbling Numbers}\label{PebNum}
Let $G=(V,E)$ be a connected graph with $n=n(G)$ vertices
$V=\{v_1,\ldots,v_n\}$ and having edge set $E$.
A {\it configuration} $C$ of $t$ pebbles on $G$ is an assignment of $t$
indistinguishable pebbles to the vertices of $G$.
The notation $C(v)=x$ means that $x$ pebbles have been assigned to the
vertex $v$.
(Notation such as $C(a,b,c)=x$ means that $x$ pebbles have been assigned
to each of the vertices $a,b,c$.)
The parameter $t$ is known as the {\it size} of $C$, which is also denoted
by $|C|$.
A {\it pebbling step} from vertex $u$ to vertex $v$ involves the removal
of two pebbles from $u$ and the placement of one of them onto $v$.
A configuration is $r$-{\it solvable} if it is possible to place a
pebble on the {\it root} vertex $r$ via a (possibly empty) sequence
of pebbling steps.
A configuration is {\it solvable} if it is $r$-solvable for all choices
of a root $r$.
The {\it pebbling number} of a graph $G$, denoted $\pi(G)$, is the least
number $t$ for which every configuration of size $t$ is solvable.
One should read \cite{Hur2} for the history and main results of the
theory of pebbling in graphs.

Clearly $\pi(G)\ge n(G)$ for every $G$, for the configuration which places
no pebble on $r$ and one pebble on all other vertices is not $r$-solvable.
The authors of \cite{CHH} defined $G$ to be {\it Class 0} when $\pi(G)=n(G)$.
Examples of Class 0 graphs include cliques (via the Pigeonhole Principle)
and cubes (see \cite{Chu}), as well as the 5-cycle, the Petersen graph,
and many others.
They proved that all 3-connected, diameter 2 graphs are Class 0, and
conjectured that fixed diameter graphs with large enough connectivity
are also Class 0.
This conjecture was proved in \cite{CHKT}, where the result was used to
prove that the random graph in which each edge appears independently with
probability $p\gg (n\lg n)^{1/d}/n$ for some $d>0$ is almost surely Class 0.

Let $\bB(m)$ be the set of all connected bipartite graphs with $m$
vertices in each part.
It was proven in \cite{Hur1} that every regular graph in $\bB(m)$ having
degree at least $\lc 2m/3\rc+1$ is Class 0.
In this paper we derive a more general and stronger result for large $m$.
Let $b=b(m)$ be the minimum number so that every $G\in\bB(m)$ having
minimum degree at least $b$ is Class 0.

\begin{theorem}\label{BT}
For all $m\ge 336$, $b(m)=\lf m/2\rf+1$.
\end{theorem}

We prove this in Section \ref{BP}.
We also derive an analogous result for connected graphs.
Let $\bG(n)$ be the set of all connected graphs on $n$ vertices,
and let $g(n)$ be the minimum number $g$ so that every $G\in\bG(n)$
having minimum degree at least $g$ is Class 0.
We prove the following in Section \ref{GP}.

\begin{theorem}\label{GT}
For all $n\ge 9$, $g(n)=\lf n/2\rf$.
\end{theorem}

%
%
\subsection{Thresholds}\label{Thresh}
We next consider a randomized version of pebbling, introduced in
\cite{Cla}, in which we consider the probability space of all
configurations of $t$ pebbles, each equally likely.
The pebbling number is the minimum $t$ for which the probability
that a configuration is solvable equals 1.
Now we wish to find $t$ so that this probability is nearly 1.
To be more precise, let us introduce some notation.

For two functions $f=f(n)$ and $g=g(n)$ we say that $f\ll g$ ($g\gg f$)
if $f/g\ra 0$ as $n\ra\inf$.
We set $o(g)=\{f\ |\ f\ll g\}$ and $\w(f)=\{g\ |\ g\gg f\}$.
We also write $f\sim g$ whenever $f/g\ra 1$ as $n\ra\inf$.
Further, we set
$O(f)=\{g\ |\ {\rm for\ some\ } c,k>0, g<cf {\rm\ for\ all\ } n>k\}$,
and similarly
$\W(g)=\{f\ |\ {\rm for\ some\ } c,k>0, f>cg {\rm\ for\ all\ } n>k\}$.
Finally we define $\Th(f)=O(f)\cap\W(f)$.

We consider sequences $\G=(G_1,\ldots,G_n,\ldots)$ of graphs for which
the number of vertices increases with $n$ (e.g. $G_n$ has $n$ vertices).
For a function $t=t(n)$ we denote by $\Pr_t(n)$ the probability that
a randomly chosen configuration of $t$ pebbles on $G_n$ is solvable.
A function $\t=\t(n)$ is a {\it pebbling threshold} for $\G$ if
$\Pr_t(n)\ra 1$ for all $t\gg\t$ and $\Pr_t(n)\ra 0$ for all $t\ll\t$.
We denote by $\t(\G)$ the set of all pebbling thresholds for $\G$.
It is not immediately evident that every graph sequence has such a threshold,
but it is proven so in \cite{BBCH}.
The first threshold result, from \cite{Cla}, established that 
$\t(\K)=\Th(n^{1/2})$, where $\K$ is the sequence of complete graphs.
More results on the pebbling thresholds of paths, cubes, and other sequences
appear in \cite{BBCH,CEHK,CW,GJSW}.
For instance, it is known that if $t\in\t(\G)$ for some graph sequence $\G$,
then $t\in\W(n^{1/2})$.
\vspace{0.2 in}

For our purposes let us define $\bG(n,\d)$ to be the set of all connected
graphs on $n$ vertices having minimum degree at least $\d=\d(n)$.
Let $\G_\d=(G_1,\ldots,G_n,\ldots)$ denote any sequence of graphs with
each $G_n\in\bG(n,\d)$.
In Section \ref{DP2} we prove the following theorem.

\begin{theorem}\label{DT2}
For every function $\sqrt{n} \ll \d=\d(n)\le n-1$, $\t(\G_\d)\subseteq O(n^{3/2}/\d)$.
In particular, if in addition $\d\in\W(n)$ then $\t(\G_\d)=\Th(n^{1/2})$.
\end{theorem}

%
%
\section{Proofs}\label{Proofs}

%
%
\subsection{Theorem \ref{BT}}\label{BP}
\noindent {\bf Lower bound.}
First we give a proof of the lower bound, that $b(m)\ge\lf m/2\rf +1$
for all $m\ge 7$.

We define, for each $m$, the bipartite graph $B_m$ as follows.
Let $|L|=|R|=m$ with $L=L_1\cup L_2$ and $R=R_1\cup R_2$ so that
$|L_1|=|R_1|=\lc m/2\rc$ and $|L_2|=|R_2|=\lf m/2\rf$.
Let the induced subgraphs on $L_1\cup R_1$ and on $L_2\cup R_2$ each be
complete bipartite with one missing edge.
Let the two missing edges be $xy$, with $x\in L_1$ and $y\in R_1$,
and $wz$, with $w\in L_2$ and $z\in R_2$.
Finally include the two edges $wy$ and $xz$.
Note that for $\d=\lf m/2\rf$ the graph $B_m$ has minimum degree $\d$,
and is $\d$-regular when $m$ is even.

Now we define a configuration $C$ of size $n=n(B_m)=2m$, and show that it is
unsolvable when $m\ge 7$.
We choose the root $r\in L_2-\{w\}$ and define $C(r,w,x,y,z)=0$.
We find $a,b,c\in L_1-\{x\}$ and define $C(a,b)=3$ and $C(c)=2$.
Finally we define $C(v)=1$ for all other vertices $v$.
Clearly, $|C|=n$.

In order that $C$ is $r$-solvable one must be able to move 2 pebbles onto
either $w$ or $z$, and consequently 4 pebbles onto either $x$ or $y$.
It is not difficult to see that both cases are impossible, since at most
three pebbles can be put in motion via pebbling steps from $a,b$, and $c$.
Hence, for all $m\ge 7$, $B_m$ is not Class 0 and so $b(m)\ge\lf m/2\rf +1$.
\hfill\pf
\vspace{0.2 in}

\noindent {\bf Upper bound.}
Second we give a proof of the upper bound, that $b(m)\le\lf m/2\rf+1$
for all $m\ge 336$.

Let $B\in\bB(m)$ have bipartition $L,R$ and 
minimum degree at least $\lf m/2\rf +1$, where $|L|=|R|=m$.
Choose any configuration $C$ of size $n=n(B)=2m$ and let $r$ be any chosen
root (which we may assume lies in $L$).
We assume that $C$ is $r$-unsolvable and derive a contradiction.

We will make use of the following two observations about $B$.
First, every pair of vertices in the same part has a common neighbor.
Second, from this it is clear that the diameter of $B$ is at most 3.
We will derive a contradiction by accumulating 8 pebbles on some vertex,
from which we can obviously pebble to $r$.

Denote the neighborhood of a vertex $v$ by $N(v)$, and the union of
neighborhoods of a set $S$ of vertices by $N(S)$.
Then we must have $C(r)=0$ and $C(v)\le 1$ for all $v\in N(r)$.
We also know that $C(v)\le 3$ for all $v\in N(N(r))=L$.
Let $Z=\{v\ |\ C(v)=0\}$, $U=\{u\ |\ C(u)=1\}$, and $H=V-Z-U$.
We let $Z_L=Z\cap L$, with $Z_R$, $U_L$, $U_R$, $H_L$, and $H_R$
defined analogously.
\vspace{0.2 in}

\noindent
{\bf Claim.}
$|Z|>m/2$.
\vspace{0.2 in}

\noindent
{\it Proof.}
The claim is trivial if $N(r)\subseteq Z$ so we assume otherwise and
pick some $r\"\in N(r)\cap U$.
Now we know that we cannot move another pebble to $r\"$.
We note that $H$ is nonempty because $Z$ is nonempty.
Also we note that if there is some vertex $v$ with $C(v)\ge 4$,
then we can move a pebble to either $r$ or $r\"$.
Hence we assume $C(v)\le 3$ for all $v$.
Moreover, we note that it must be impossible to ever put 4 pebbles
on any vertex.

Suppose that there is a vertex $t$ having $C(t)=3$; without loss of
generality we assume that $t\in L$ (if $t\in R$, we think of $r\"$
as our new root and argue similarly).
Since $N(r)\cap N(t)\subseteq Z$, we know that $|Z|\ge 2$.
Moreover, $N(t)\subseteq Z$ implies $|Z|>m/2$,
so we assume otherwise and pick $t\"\in N(t)\cap (U\cup H)$.
Then $N(t\")\cap N(r\")\subseteq Z$, and so $|Z|\ge 3$,
which implies that $|H|\ge 2$.

If there is a vertex $s\in H_L-\{t\}$ then we can argue as follows.
Let $X=N(r)$, $T=N(t)$, and $S=N(s)$.
Of course, $(X\cap S)\cup (X\cap T)\cup (S\cap T)\subseteq Z$
($S\cap T\subseteq Z$ since otherwise we could place 4 pebbles on $t$),
and $X\cap S\cap T=\mt$.
Therefore we have that
\begin{center}
\begin{tabular}{rcl}
$m$&$\ge$&$|X\cup S\cup T|$\\
   &$=$&$|X|+|S|+|T|-|X\cap S|-|X\cap T|-|S\cap T|$\\
   &$>$&$3m/2-|Z|$,\\
\end{tabular}
\end{center}
which implies that $|Z|>m/2$.
If instead there is a vertex $s\"\in H_R$ then we know that,
either $|Z|>m/2$ because $N(s\")\subseteq Z$,
or there is some $s\in N(s\")\cap (U\cup H)-\{t\}$.
In the latter case we move a pebble from $s\"$ to $s$ and argue as above.
Henceforth we may assume that $C(v)\le 2$ for all $v$.

Consequently the equality
$$|Z|+|U|+|H| = n(B) = |C| = |U|+2|H|$$
tells us that $|H|=|Z|$, and so $|H|\ge 3$.
Therefore, again without loss of generality, $|H_L|\ge 2$,
say $\{s,t\}\subseteq H_L$.
If $S\cap T\subseteq Z$ then we may copy the above argument that
$m\ge |X\cup S\cup T|$ implies $|Z|>m/2$.
Otherwise we may move a pebble from $s$ through $S\cap T$ to $t$,
find $p\in H-\{s,t\}$ and use the original argument for the case that
$C(t)=3$.
This completes the proof of the Claim.
\hfill\cpf

Now we can use the relations 
$$|Z|+|U|+|H| = n(B) = |C| = |U|+\sum_{v\in H}C(v)$$
to see that
$${m\over 2} < |Z| = \sum_{v\in H}C(v) -|H| \le 6|H|\ ,$$
so that $|H|>m/12$.  From this we can assume, without loss of generality, 
that $|H_L|>m/24$,
so that the number of edges with one end in $H_L$ is more than $(m/24)(m/2)$.
Since $m\ge 336$ there must be some $x\in R$ having at least 8 neighbors
in $H_L$, so that we can put 8 pebbles on $x$, a contradiction.

This contradiction proves that $C$ is $r$-solvable.
\hfill\pf

%
%
\subsection{Theorem \ref{GT}}\label{GP}
\noindent {\bf Lower bound.}
First we give a proof of the lower bound, that $g(n)\ge\lf n/2\rf$
for all $n\ge 9$.

We define, for each $n$, the graph $G_n$ as follows.
Let the vertex set $V=L\cup R$, with $|L|=\lc n/2\rc$ and $|R|=\lf n/2\rf$.
Let the induced subgraphs on $L$ and on $R$ each be
complete with one missing edge.
Suppose the edge $xy$ is missing from the subgraph on $L$,
and the edge $wz$ is missing from the subgraph on $R$.
Finally include the two edges $wy$ and $xz$.
Note that for $\d=\lf n/2\rf-1$ the graph $G_n$ has minimum degree $\d$,
and is $\d$-regular when $n$ is even.

Next we define a configuration $C$ of size $n=n(G_n)$,
and show that it is unsolvable when $n\ge 9$.
We choose the root $r\in R-\{w,z\}$ and define $C(r,w,x,y,z)=0$.
We find $a,b,c\in L-\{x,y\}$ and define $C(a,b)=3$ and $C(c)=2$.
Finally we define $C(v)=1$ for all other vertices $v$.
Clearly, $|C|=n$.

In order that $C$ is $r$-solvable one must be able to move 2 pebbles onto
either $w$ or $z$, and consequently 4 pebbles onto either $x$ or $y$.
It is not difficult to see that both cases are impossible.
Hence for all $n\ge 9$, $G_n$ is not Class 0, and so $b(n)\ge\lf n/2\rf$.
\hfill\pf
\vspace{0.2 in}

\noindent {\bf Upper bound.}
Second we give a proof of the upper bound, that $g(n)\le\lf n/2\rf$
for all $n\ge 6$.

\begin{figure}
\centerline{\hbox{\psfig{figure=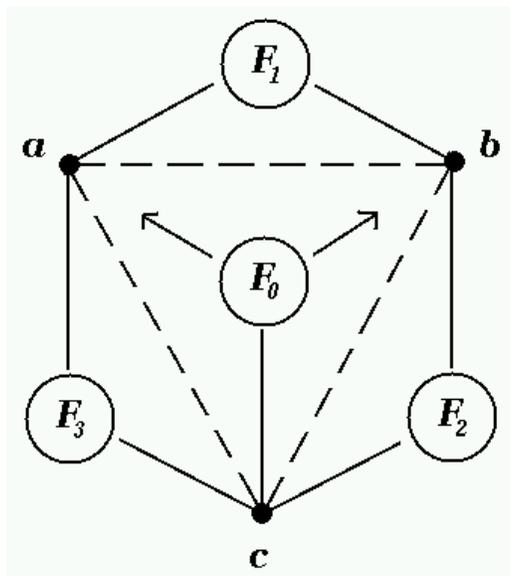,height=3.0in}}}
\caption{a schematic diagram of Class 1 graphs of diameter two}
\label{Class1}
\end{figure}

Let $G$ be graph with minimum degree $\lf n/2\rf$.
We suppose that $G$ is not Class 0 and derive a contradiction.
Because complete graphs are Class 0, $G$ has diameter at least 2,
and because every pair of vertices of $G$ has a common neighbor,
the diameter of $G$ is exactly 2.
It is proven in \cite{PSV} that every graph $G$ of diameter two has pebbling
number $n(G)$ or $n(G)+1$ ({\it Class 1}).
In \cite{CHH} we find the following characterization of 
Class 1 graphs of diameter two (see Figure \ref{Class1}).

In the figure, $F_0$ is any (possibly empty) graph, $F_1$ is any
nonempty graph, and $F_2$ and $F_3$ are any nonempty connected graphs.
The solid lines indicate an edge from the given vertex to every vertex
in the corresponding set.
At least two of the three dashed lines must be present, and the arrows
indicate that every vertex in $F_0$ must have at least one edge to $\{a,b\}$.

Suppose that $G$ is labeled as in Figure \ref{Class1}.
Choose any vertices $p\in F_3, q\in F_2, r\in F_1$, and recall that
$N(x)$ denotes the neighborhood of a vertex $x$.
Since $\{p,q\}$ is not an edge, $|N(p)\cap N(q)|=1$ and
$r\notin N(p)\cup N(q)$, we have $deg(p)+deg(q)\le n-2$, and so at least one
of $p$ or $q$ has degree smaller than $\lf n/2\rf$, a contradiction.
Hence $G$ is Class 0.
\hfill\pf

In fact it is easy to show by induction that every diameter two Class 1
graph has minimum degree at most $\lf n/3\rf$.

%
%
\subsection{Theorem \ref{DT2}}\label{DP2}
In this section we prove Theorem \ref{DT2}.
The proof is divided into two steps.
First, we will show that it is possible to partition the vertices of 
$G_n=(V,E)$ into $O(\frac{n}{\d})$ subgraphs of diameter at most two. 
Second, we prove that if a distribution has enough pebbles then there
will be in every subgraph ``many'' vertices with two pebbles each.

We call a partition $V_1,\dots,V_l,W$ a $q$-{\it star partition} of $V$ if
\begin{enumerate}
\item for every $1 \leq i \leq l$, $V_i$ contains a star on 
at least $q$ vertices and
\item every vertex of $W$ has a neighbor in $V_i$ for some $1 \leq i \leq l$.
\end{enumerate}
The following procedure constructs a $(\d+1)$-star partition 
$V_1, \dots , V_l, W$ of $V$ with $l=O(n/\d)$.
Select $v \in V$ arbitrarily and let $V_1= N(v) \cup \{v\}$. 
For a general step, suppose $V_1, \dots , V_k$ have been selected and let
$U = V \setminus \bigcup_{i=1}^k V_i$.
Either every vertex from $U$ has a neighbor in $V_i$ for some
$1 \leq i \leq k$, in which case we stop the process with $l=k$ and $W=U$,
or there is a vertex $u \in U$ such that $N(u) \subseteq U$.
In the latter case we set $V_{i+1}=N(u)\cup \{u\}$ and continue the process.

Assume that $V_i=\{v_{i0},v_{i1},\dots, v_{ik_i}\}$ with $k_i \geq \d$
and let $C$ be a pebbling distribution with $t = \w n^{3/2}/\d$
pebbles where $\w =\w(n) \rightarrow \infty$ is such that $t \leq n-1$.
Consider the random variable $X_{ij}$, which is equal to one if $C(v_{ij})=2$
and zero otherwise, and let $X_i =\sum_{j=0}^{k_i}X_{ij}$.
Then
$$ \E[X_{ij}] = \frac{{t+n-4 \choose t-2}}{{t+n-1 \choose t}}\ ;$$
using the assumptions about $t$ and $n$ it is easy to
check that $\E[X_{ij}X_{ik}] \leq \E[X_{ij}]\E[X_{ik}]$ for $ j \neq k$.
Consequently
$$var[X_i] = \E[X_i^2] - \E[X_i]^2 \leq \E[X_i]\ ,$$
and by Chebyshev's Inequality
$$\Pr\[X_i < \frac{\E[X_i]}{2}\]
	\leq  \Pr\[|X_i - \E[X_i]| \geq \frac{\E[X_i]}{2}\]
	\leq \frac{4}{\E[X_i]}\ .$$
Thus the probability that there exists an $1\leq i \leq l$ such that
$X_i < \E[X_i]/2$ is $O(n/\d\E[X_i])$.
But, with our choice of $t$ and $n$,
$$\E[X_i] = \W\(\frac{\d t^2}{n^2}\)
	= \W\(\frac{\w^2 n}{\d}\)\ .$$
Thus, with probability tending to one, for all $i$, we have
$X_i \geq \E[X_i]/2$, which is at least 8 for large enough $n$.
Therefore, with probability tending to one, 4 pebbles can be accumulated
on the center of every star, and since every vertex is within distance two
of some center, it is possible to move a pebble to any given root vertex.

%
%

%
%
\bibliographystyle{plain}
%

%

%
%
\end{document}